\documentclass{article}

\usepackage{arxiv}

\usepackage[utf8]{inputenc} 
\usepackage[T1]{fontenc}    
\usepackage{hyperref}       
\usepackage{url}            
\usepackage{booktabs}       
\usepackage{amsfonts}       
\usepackage{nicefrac}       
\usepackage{microtype}      
\usepackage{lipsum}
\usepackage{graphicx}
\graphicspath{ {./images/} }
\usepackage{siunitx}
\usepackage[utf8]{inputenc}
\usepackage{placeins}
\usepackage{times}
\usepackage{multirow}
\usepackage{amsmath}
\usepackage{graphicx}
\usepackage{tabularx}
\usepackage{lineno}
\usepackage[export]{adjustbox}
\usepackage{amsmath}
\usepackage{epstopdf}
\usepackage{xcolor}
\definecolor{amber}{rgb}{1.0, 0.75, 0.0}
\definecolor{brickred}{rgb}{0.7960, 0.2550, 0.3290}
\usepackage{placeins}
\usepackage{times}
\usepackage{multirow}
\usepackage{amsmath,amssymb}

\newcommand{\subfigimgthree}[3][,]{%
  \setbox1=\hbox{\includegraphics[#1]{#3}}
  \leavevmode\rlap{\usebox1}
  \rlap{\hspace*{-10pt}\raisebox{\dimexpr\ht1-1\baselineskip}{#2}}
  \phantom{\usebox1}
}

\usepackage{tikz}

\DeclareMathOperator*{\argmax}{\arg\!\max}
\usepackage{subcaption}
\usepackage{mathtools}
\usepackage{bm}
\usepackage[utf8]{inputenc}
\usepackage{booktabs}

\newcommand{\citet}[1]{\citeuthor{#1} \shortcite{#1}}

\usepackage[utf8]{inputenc}
\usepackage{placeins}
\usepackage{times}
\usepackage{multirow}
\usepackage{amsmath}
\usepackage{graphicx}
\usepackage{algorithmic}
\usepackage{algorithm}
\usepackage{lineno}
\usepackage[export]{adjustbox}
\usepackage{amsmath}
\usepackage{epstopdf}
\usepackage{xcolor}
\usepackage{bm}
\definecolor{amber}{rgb}{1.0, 0.75, 0.0}
\definecolor{brickred}{rgb}{0.7960, 0.2550, 0.3290}
\definecolor{seagreen}{rgb}{0.18, 0.55, 0.34}
\DeclareMathOperator*{\argmin}{arg\,min}

\usepackage{placeins}
\usepackage{times}
\usepackage{multirow}
\usepackage{amsmath}
\usepackage{graphicx}
\usepackage{lineno}
\usepackage[export]{adjustbox}
\usepackage[resetlabels,labeled]{multibib}
\newcites{S}{Supplemental References}

\title{Collective wind farm operation based on a predictive model increases utility-scale energy production}

\author{
  Michael F. Howland \\
  Civil and Environmental Engineering \\ Massachusetts Institute of Technology \\ Cambridge, MA 02139, USA \\
  \texttt{mhowland@mit.edu} \\
  \And
  Jes{\'u}s Bas Quesada, Juan Jos{\'e} Pena Mart{\'i}nez, Felipe Palou Larra\~{n}aga \\
  Siemens Gamesa Renewable Energy Innovation \& Technology \\
  Sarriguren, Navarra, Spain, 31621 \\
  \And
  Neeraj Yadav, Jasvipul S. Chawla, Varun Sivaram \\
  ReNew Power Private Limited \\
  Gurugram-122009, Haryana, India \\
  \And 
  John O. Dabiri \\
  California Institute of Technology \\
  Pasadena, CA 91125, USA \\
}

\begin{document}
\maketitle
\begin{abstract}
Wind turbines located in wind farms are operated to maximize only their own power production. 
Individual operation results in wake losses that reduce farm energy.
In this study, we operate a wind turbine array collectively to maximize total array production through wake steering.
The selection of the farm control strategy relies on the optimization of computationally efficient flow models.
We develop a physics-based, data-assisted flow control model to predict the optimal control strategy.
In contrast to previous studies, we first design and implement a multi-month field experiment at a utility-scale wind farm to validate the model over a range of control strategies, most of which are suboptimal.
The flow control model is able to predict the optimal yaw misalignment angles for the array within $\pm 5^\circ$ for most wind directions ($11$--$32\%$ power gains).
Using the validated model, we design a control protocol which increases the energy production of the farm in a second multi-month experiment by $2.7\%$ and $1.0\%$, for the wind directions of interest and for wind speeds between $6$ and $8~\mathrm{m/s}$ and all wind speeds, respectively.
The developed and validated predictive model can enable a wider adoption of collective wind farm operation.
\end{abstract}

\section*{Introduction}

With consensus that carbon-intensive energy generation has contributed to global warming \cite{ipcc_ar6_2021}, the decarbonization of electricity production is of paramount importance \cite{ipcc_1p5_deg}.
Renewable energy must produce more than $60\%$ of the primary energy supply by 2050 to achieve the $2^\circ$C warming target set by the Paris Agreement \cite{gielen2019role}.
Just $14\%$ of worldwide energy production in 2015 was from renewables \cite{gielen2019role}.
Meanwhile, the United States has set a goal of $100\%$ carbon-free electricity generation by 2035 \cite{bide2021renewable}, with similar international goals elsewhere \cite{renn2016coal}.
Further, approaches to maximize renewable energy production in emerging economies are critical to address climate change \cite{asif2007energy}.
Improvements in the efficiency of wind generation will enable a more rapid and lower cost transition to a decarbonized energy system \cite{veers2019grand}.

While individual horizontal axis wind turbines are approaching theoretical peak efficiency \cite{wiser2015wind}, wind farms exhibit losses from turbine interactions.
Utility-scale wind farms lose $10$--$20\%$ of their energy production per year due to wake interactions between turbines \cite{barthelmie2009modelling}.
Individual turbines generate turbulent, energy deficit wake regions downwind \cite{stevens2017flow}.
Wind turbines are placed in close proximity to decrease the levelized cost of energy (LCOE) for the collective farm by reducing the capital costs \cite{eberle2019nrel}.
The result is that modern turbines are spaced $6-10$ rotor diameters apart in onshore wind farms \cite{stevens2017combining} which results in significant wake interactions \cite{meyers2012optimal}.
All utility-scale wind turbines are operated to maximize their individual power production \cite{pao2009tutorial}.
Such control inherently neglects wake interactions between neighboring turbines.

Initial field experiments of collective operation at utility-scale wind farms implemented operational strategies which are based on the optimization of a flow control model \cite{howland2019wind, fleming2019initial, doekemeijer2021field}. 
These studies have demonstrated that collective operation can increase power production for wind conditions which result in high wake losses, compared to individual control.
However, to achieve the maximum power production of the wind farm, the flow control models must reliably predict the optimal control strategy.
It has not yet been demonstrated, using field data, that flow control models are able to predict the optimal control strategy for utility-scale wind farms, since this requires sustained operation in suboptimal strategies.

\begin{figure}
    \includegraphics[width=1.0\linewidth]{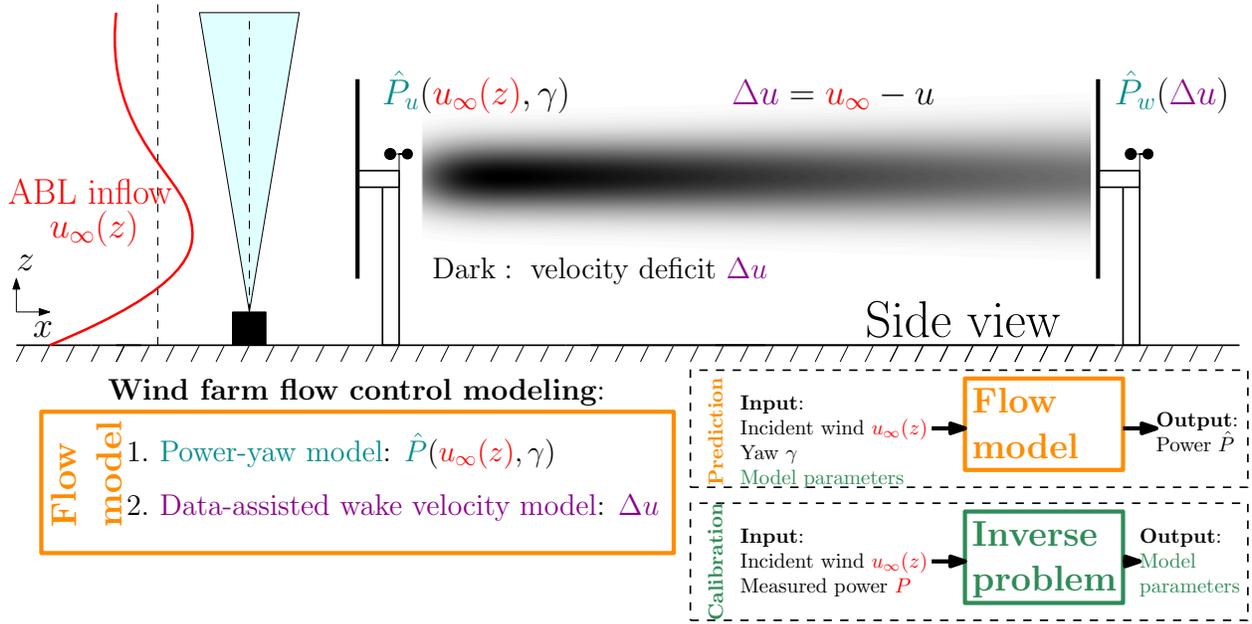}
    \caption{
    Schematic of the predictive flow control model.
    The atmospheric boundary layer (ABL) velocity profile is illustrated in red. 
    Incident ABL winds $u_\infty(z)$ are measured in the field experiments using a LiDAR (black box and blue cone). 
    Each turbine is equipped with cup and sonic anemometers (black circles), and generates a wake region (dark shaded area).
    To predict the effect of a control strategy on the power of the collective wind farm, we model the power production of upwind turbines operating in freestream conditions $\hat{P}_u$ and the waked turbines $\hat{P}_w$.
    The flow control model proposed in this study is the combination of a power-yaw model $\hat{P}(u_\infty(z),\gamma)$, which predicts the power production of a yawed turbine based on the incident wind $u_\infty(z)$ and the yaw misalignment $\gamma$, and a data-assisted wake model, which predicts the wake velocity deficit $\Delta u$.
    The wake model parameters are calibrated using wind farm data for which the turbines are operating in baseline, yaw aligned conditions $\bm{\gamma}=0$.
    The wake model is then used to predict the farm power given a yaw control strategy $\bm{\gamma}\neq0$.
    }
    \label{fig:model}
\end{figure}

In this study, we develop a methodology for operating wind turbines collectively to maximize wind farm energy production based on a new, predictive wind farm flow control model (described in {\bf Methods}).
Importantly, we conduct our experiments in configurations predicted to result in optimal performance as well as suboptimal regimes. 
This enables direct validation of the flow control model predictions.

We develop and implement a two phase field experiment at a utility-scale wind farm in India.
In the first phase, we modify the operation of a wind turbine over a broad range of control strategies, most of which are suboptimal, to identify the optimal strategy.
We demonstrate, for the first time, that our flow control model is able to predict the true optimal control strategy for a utility-scale wind farm.
We observe power production gains between $11$--$32\%$ using collective operation for wind speeds between $6$ and $8~\mathrm{m/s}$, compared to individual control.
In the second phase, we use the validated flow control model to design an optimal control strategy for the farm.
We demonstrate a $1.0\% \pm 0.5\%$ energy increase for the utility-scale wind farm for the wind directions of interest, compared to standard individual control.

\section*{Collective wind farm control background and motivation}

Individual wind turbine operation attempts to minimize the yaw angle of misalignment between the incident wind direction and the turbine rotor orientation \cite{fleming2014field}.
We consider collective wind farm operation through wake steering control, wherein certain wind turbines in the wind farm are intentionally misaligned in yaw with respect to the incident wind direction.
The power production of the yaw misaligned wind turbine is generally reduced, because the wind velocity perpendicular to the rotor is reduced \cite{howland2020influence}. 
The power production for the waked turbine may be increased due to wake deflection associated with the yawed turbine \cite{fleming2015simulation}.
The goal of wake steering optimization is to select the yaw misalignment angles which maximize the power production of the wind farm, by achieving increases in the power production of downwind turbines that compensate for the loss in power of the yaw-misaligned upstream turbines.
The optimal yaw misalignment angles inherently depend on the site-specific wind farm layout \cite{howland2019wind} and incident wind conditions \cite{fleming2017field, howland2021optimal}, which vary in time.
Therefore, the yaw optimization needs to be performed for each possible state of wind conditions \cite{quick2017optimization}, which is high-dimensional.

Since computational fluid dynamics (CFD) simulations remain intractable for such optimization \cite{choi2012grid}, the optimization of the yaw angles is generally performed with numerically efficient wake models \cite{gebraad2016wind}.
To remain tractable for optimization, wake models neglect certain flow physics \cite{meneveau2019big} and parameterize the effects of turbulence in the wind farm and the atmospheric boundary layer (ABL) \cite{howland2021wind}.
Using a wind farm operational strategy resulting from the optimization of a wake model has demonstrated potential to increase wind farm power production in large eddy simulations (LES) of idealized ABL conditions \cite{howland2021optimal, gebraad2016wind, doekemeijer2020closed}, wind tunnel experiments \cite{bastankhah2019wind, campagnolo2020wind}, and field experiments \cite{howland2019wind, fleming2019initial, doekemeijer2021field}. 

\begin{figure}
  \hspace{0.5cm}
  \centering
  \includegraphics[width=0.75\linewidth]{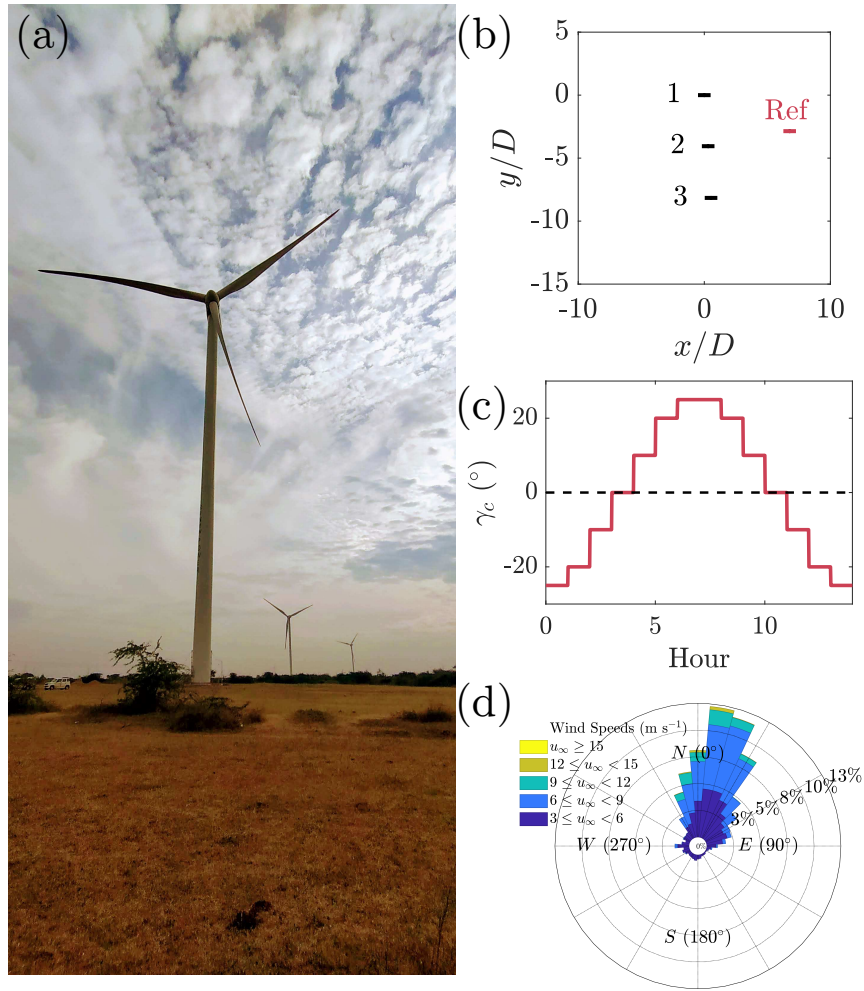}
  \caption{(a) Photo of the utility-scale wind farm of interest in this study which is located in northwest India.
  (b) Top view of the wind turbines of interest with the coordinates of the farm normalized by the wind turbine rotor diameter $D$.
  The adjacent reference turbine is denoted as `Ref'.
  (c) Commanded yaw misalignment sequence $\gamma_c$ for the fixed yaw misalignment flow control model validation experiment.
  The commanded yaw misalignments do not depend on the incident wind conditions.
  During the model validation experiment, each commanded yaw misalignment value is held fixed for one hour.
  (d) Measured wind rose during the experimental period as recorded by the reference wind turbine.
  }
    \label{fig:setup}
\end{figure}

In this study, we develop a new predictive wind farm flow control model.
We model the power production of upwind turbines given their operational strategy and given the measured incident ABL wind profiles \cite{howland2020influence}.
We model the power production of the downwind turbines using an analytical wake model \cite{howland2019wind, bastankhah2014new, shapiro2018modelling} and leveraging data-driven parameter estimation techniques \cite{howland2020optimal}.
The modeling framework is described in Figure \ref{fig:model}.

Previous field experiments of wake steering are useful to demonstrate that implementing the yaw angles resulting from the optimization of a wake model increases power compared to standard operation \cite{howland2019wind, fleming2019initial, doekemeijer2021field}.
However, it has not previously been demonstrated that the yaw angles predicted from the optimization of wake models are the true optimal yaw angles for the utility-scale wind farm.
In this study, we design and implement a field experiment at a commercial wind farm in India where we intentionally yaw misalign a wind turbine with a rotor diameter of approximately $120$ meters using fixed yaw angles between $-25^\circ$ and $25^\circ$.
We show that the proposed flow control model (Figure \ref{fig:model}) is able to predict the optimal yaw misalignment angles for the farm to within $\pm5^\circ$ for most ($5/6$) conditions tested.
Leveraging the validated wake model, we design an optimal wake steering protocol wherein the yaw misalignment angles of the wind farm vary according to the incident wind conditions.
We then perform a second field experiment where we implement the optimal wake steering protocol.

\begin{figure*}
  \centering
  \includegraphics[width=0.9\linewidth]{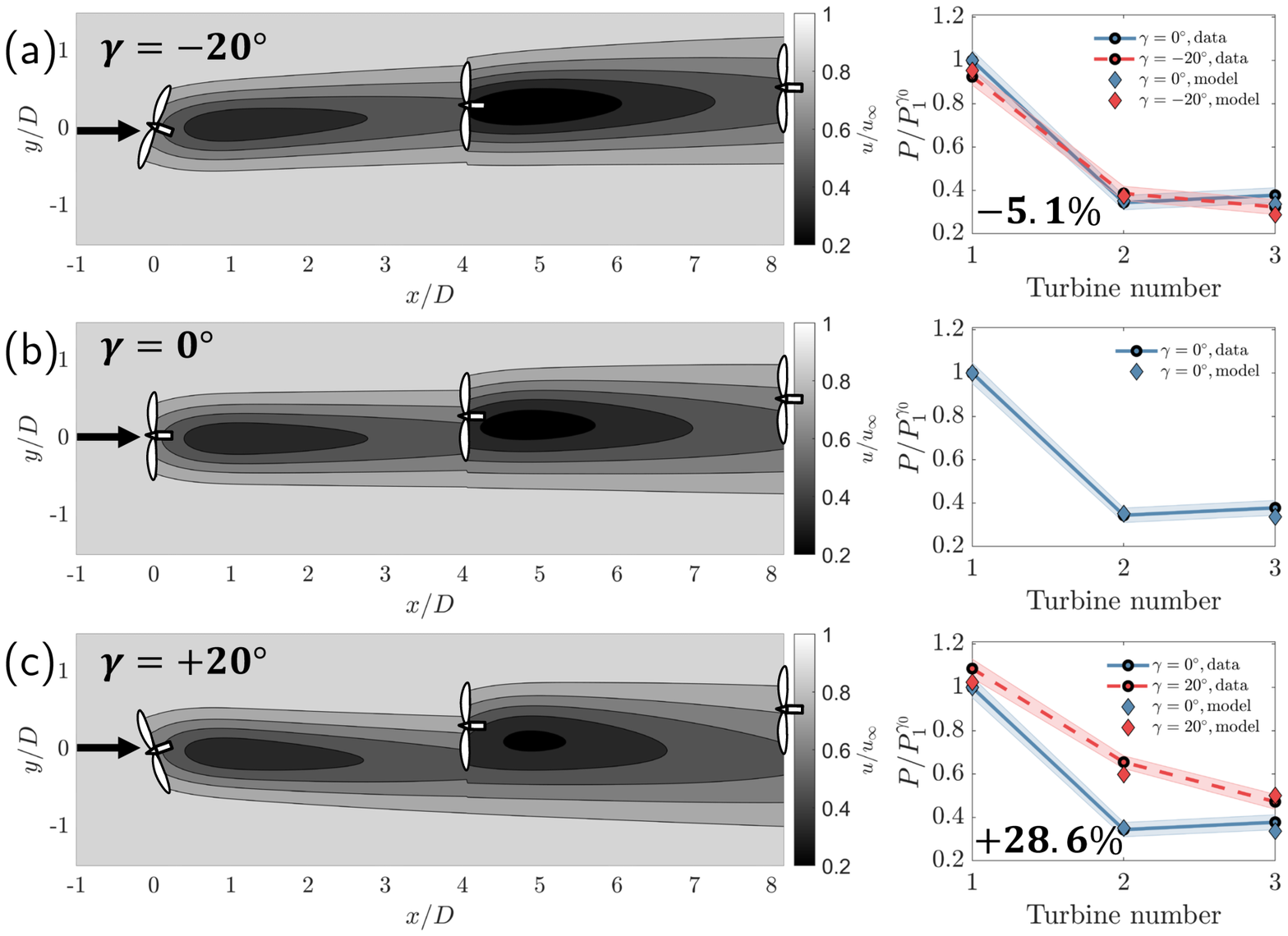}
  \caption{Results from the static yaw misalignment field experiment for turbine $1$ yaw misalignments of (a) $\gamma_1 = -20^\circ$, (b) $\gamma_1 = 0^\circ$, and (c) $\gamma_1 = 20^\circ$ for wind directions $\alpha$ from the north between $0\pm2.5^\circ$.
  The wind speed is $7 \pm 1.5$ m/s and the turbulence intensity is $5\pm2.5 \%$, as measured by the reference turbine.
  (Left) Contours of the streamwise velocity predicted by the flow control model.
  The yaw angles implemented in panels (a,c) are predicted to steer the wake region towards (a) or away (c) from the downwind turbines, depending on the yaw angle.
  Steering the wake away from the downwind turbines (c) is predicted to increase array power production (predictions in right panels).
  (Right) Comparison of the wake model predictions to the measured power production for the three turbine array.
  The powers are normalized by the power production of the leading turbine $1$ with $\gamma_1 = 0^\circ$.
  The field experiment measurements are shown as circles and the model predictions are shown as diamonds.
  The shaded region corresponds to $95\%$ confidence intervals around the mean estimated with bootstrapping.
  The blue and red correspond to yaw alignment and yaw misalignment, respectively.
  The flow control model is calibrated to the yaw aligned data $\gamma_1 = 0^\circ$.
  The flow control model is then used to predict the power production for yaw misaligned operation.
  The effect of the turbine $1$ yaw misalignment on power production of the waked turbines ($2$ and $3$) depends on the direction (sign) of the yaw.
  The field experiment measurements demonstrate a $-5.1\%$ and $+28.6\%$ change in the three turbine array mean power production compared to $\gamma_1=0^\circ$ for $\gamma_1=-20^\circ$ and $\gamma_1=+20^\circ$, respectively.
  For $\gamma_1=-20^\circ$ (a), both the data and the predictive model result in a slight decrease in the power production of turbine $1$, compared to $\gamma_1=0^\circ$. 
  For $\gamma_1=+20^\circ$ (c), both the data and the predictive model result in a slight increase in the power production of turbine $1$, compared to $\gamma_1=0^\circ$. 
  This small gain is associated with the incident wind velocity profiles in the ABL which occurred during operation with $\gamma_1=+20^\circ$; such a small gain has also been shown for certain flow conditions and yaw misalignments in previous studies \cite{doekemeijer2021field,howland2020influence}.
  }
    \label{fig:contours}
\end{figure*}

\section*{Experimental setup}

The wind farm (Figure \ref{fig:setup}(a)) consists of approximately $100$ turbines and is located in northwest India.
The turbines are approximately $100$ meters in hub height and approximately 2 MW in capacity.
We consider a subset of four turbines shown in Figure \ref{fig:setup}(b).
We focus our experiment on turbines 1, 2, and 3.
An adjacent freestream turbine is used as a reference, and is labeled as `Ref' (Figure \ref{fig:setup}(b)).
The turbines are approximately aligned for northwesterly inflow ($\alpha\approx-5^\circ$).
The prevailing wind conditions during the summer and winter are southwesterly and northeasterly, respectively, with spring and autumn as shoulder seasons with transitional winds.
We focus our experiment on the northeasterly winter winds, since the southwesterly summer winds exhibit heterogeneous upwind blockage effects, which limits the potential to perform controlled experiments.

We record data from each turbine in the form of one-minute averaged Supervisory Control and Data Acquisition (SCADA) data, which includes power production, nacelle position, blade pitch, and other relevant turbine variables.
The wind profiles depending on height are recorded in one-minute averages using a Leosphere Windcube V2.0 profiling LiDAR on site.

\subsection*{Model validation experimental setup}

We first focus on assessing the fidelity of the flow model in predicting the optimal wake steering control strategy $\bm{\gamma}^*$.
To achieve this, we implement a yaw misalignment offset series for upwind turbine 1 (Figure \ref{fig:setup}(b)) which operates in freestream conditions.
We misalign turbine 1 for yaw values between $-25^\circ$ and $25^\circ$ (Figure \ref{fig:setup}(c)).
Each yaw misalignment set-point is held fixed for one hour and does not change depending on the incident wind conditions.
Since the wind conditions vary in time, the fixed yaw set-points will often result in suboptimal operation.
We performed the model validation experiment from February 2020 until April 2020.

\begin{figure*}
  \centering
  \begin{tabular}{@{}p{0.29\linewidth}@{\quad}p{0.29\linewidth}@{\quad}p{0.29\linewidth}@{}}
    \subfigimgthree[width=\linewidth,valign=t]{(a)}{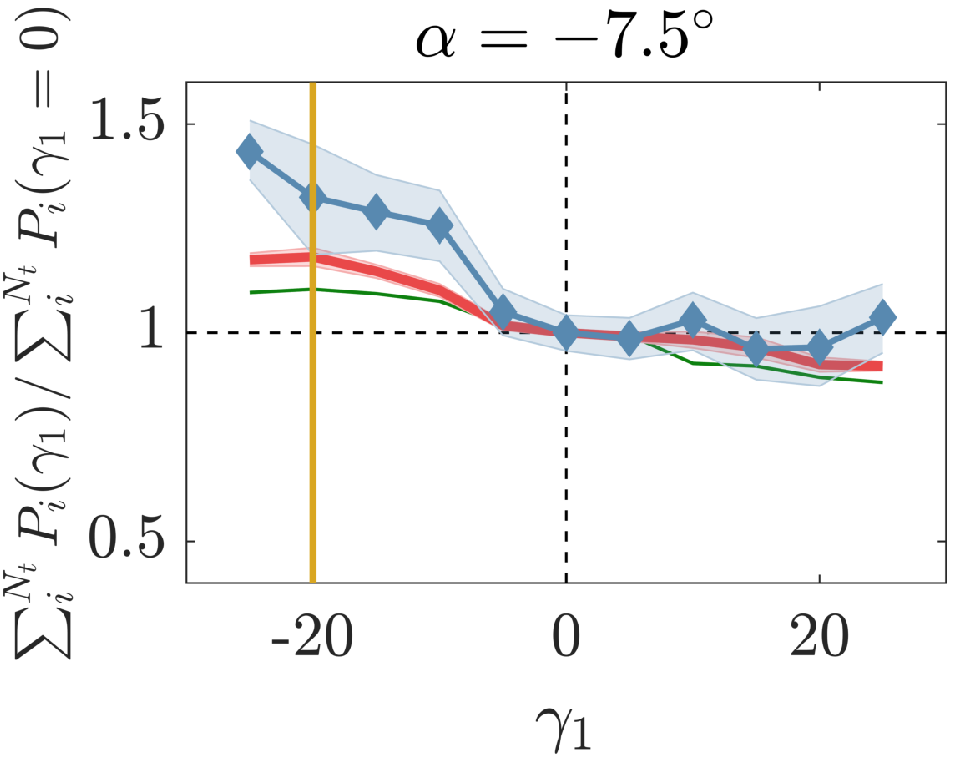} &
    \subfigimgthree[width=\linewidth,valign=t]{(b)}{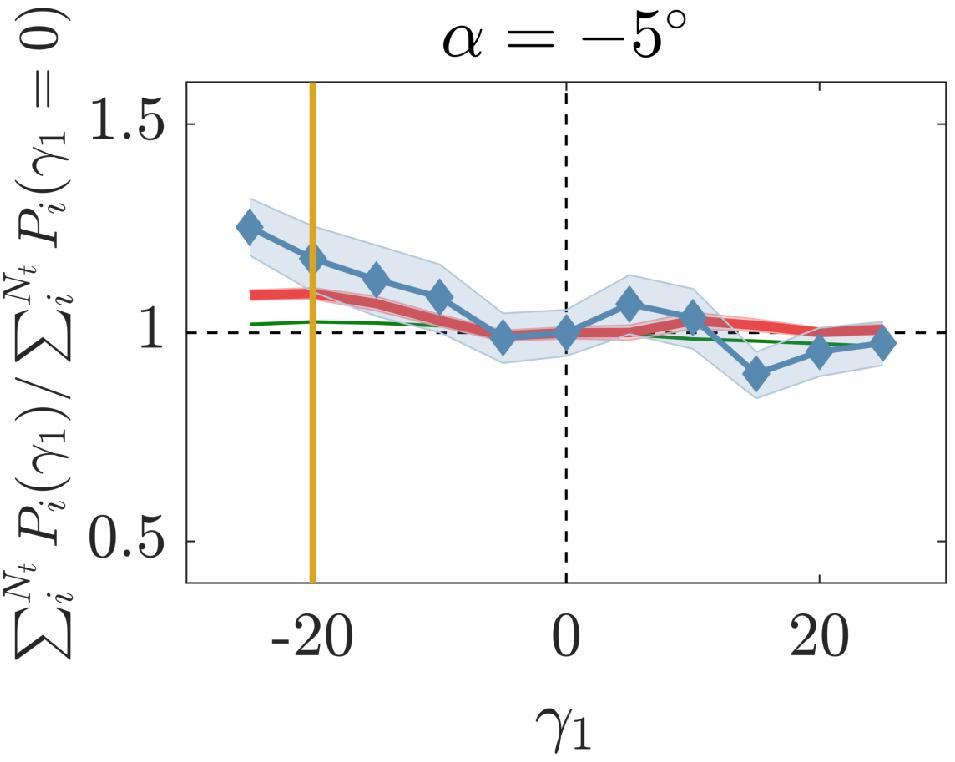} &
    \subfigimgthree[width=\linewidth,valign=t]{(c)}{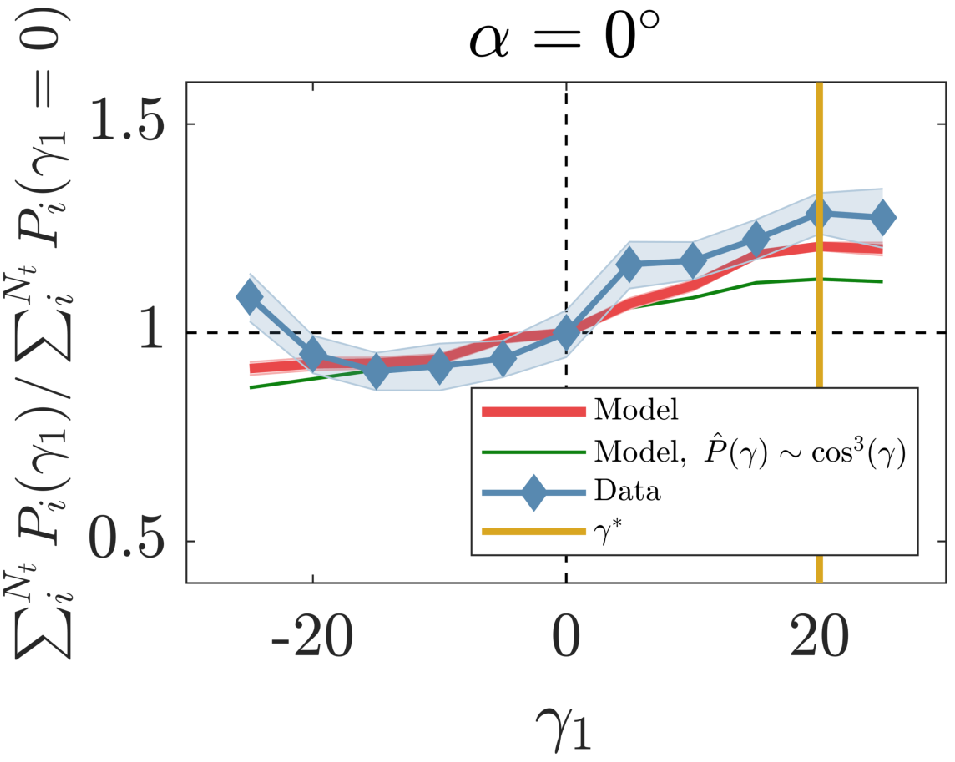} \\
    \subfigimgthree[width=\linewidth,valign=t]{(d)}{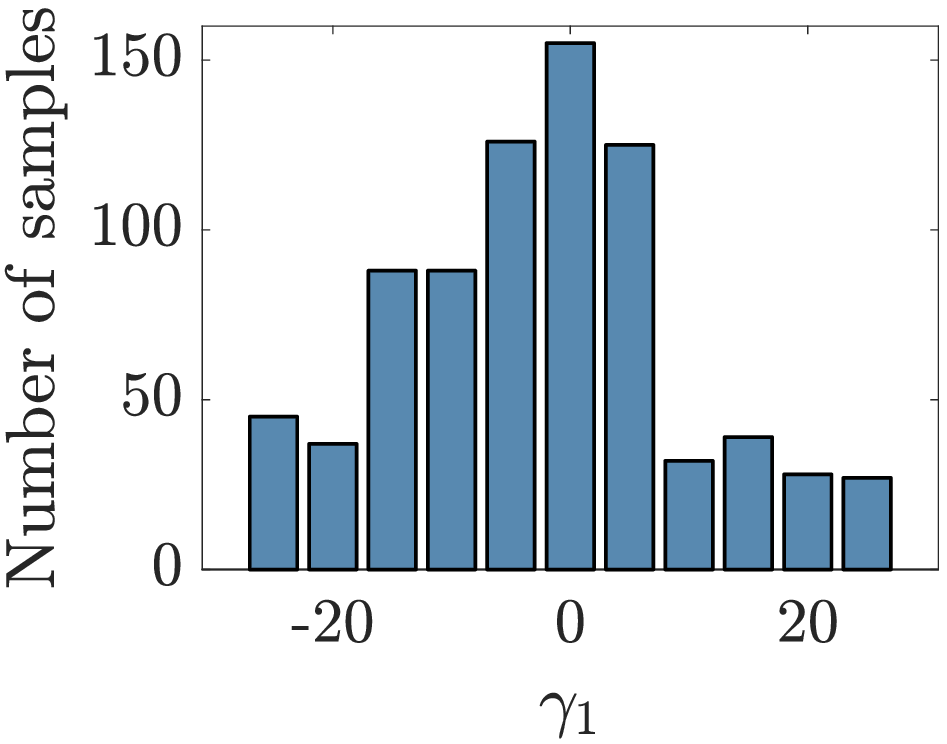} &
    \subfigimgthree[width=\linewidth,valign=t]{(e)}{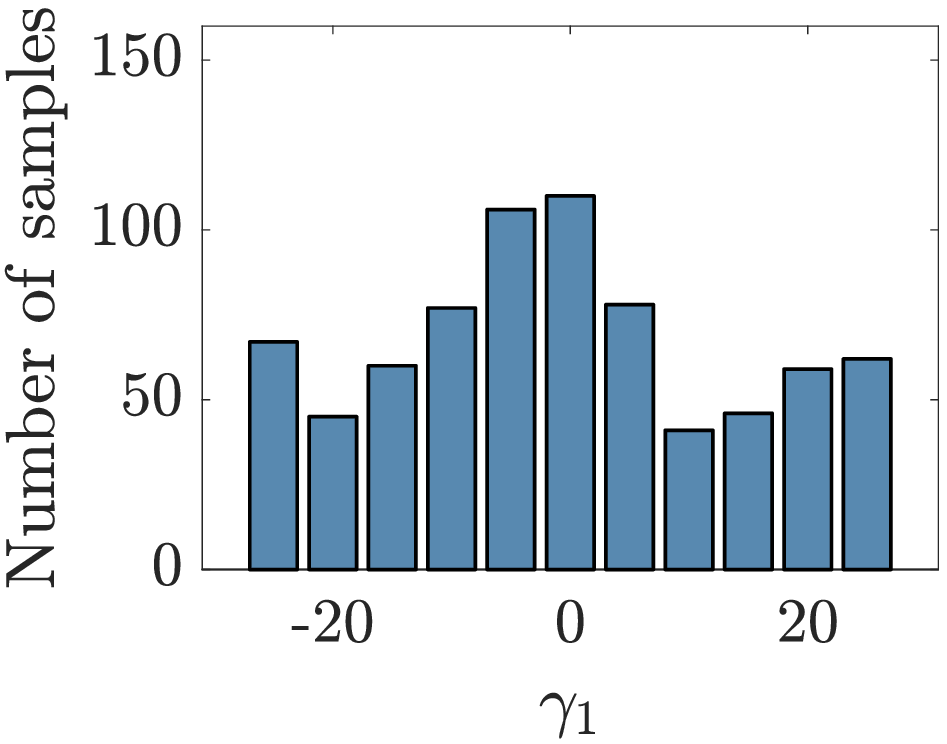} &
    \subfigimgthree[width=\linewidth,valign=t]{(f)}{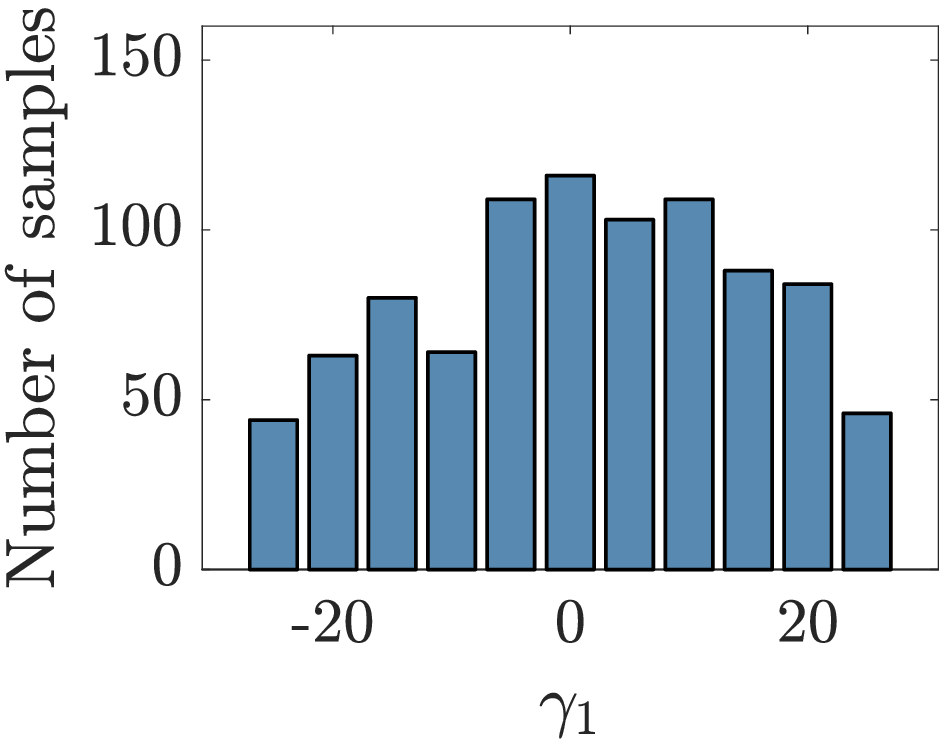}
  \end{tabular}
  \caption{Sum of the power production for turbines 1, 2, and 3 as a function of the yaw misalignment of turbine 1 from the static yaw misalignment field experiment.
  The power is shown for different inflow wind directions $\alpha$ and yaw misalignment angles $\gamma$.
  The yaw misalignment values tested were between $-25^\circ$ and $25^\circ$ (yaw values beyond $|25^\circ|$ were not considered for loading limits).
  The wind turbines are approximately aligned for northwesterly inflow ($\alpha\approx-5^\circ$).
  The power is normalized by the power produced with zero turbine 1 misalignment $\gamma_1=0^\circ$.
  The wind speed is $7 \pm 1.5$ m/s and the turbulence intensity is $5\pm2.5 \%$, as measured by the reference turbine.
  Conditional averages with $n>25$ data samples are considered.
  In blue, we show the field experiment data with $95\%$ confidence intervals from bootstrapping.
  Flow control model predictions with $95\%$ confidence intervals from bootstrapping are shown in red.
  The optimal yaw misalignment angle for turbine $1$ predicted by the flow control model is given in gold.
  Flow control model predictions assuming that the power production of yawed turbines is $\hat{P}(\gamma)\sim\cos^3 (\gamma)$ are shown in green.
  Results for other incident wind directions $\alpha$ are similar and are shown in the {\bf Supplemental Information} for brevity.
  }
    \label{fig:yaw_offset_1_2_3}
\end{figure*}

\subsection*{Farm energy maximization experimental setup}

Second, we focus on increasing the wind farm energy production through wake steering control.
The optimal yaw misalignment angles for the array depend on the incident wind conditions.
Contrary to the first experiment, we dynamically adapt the yaw misalignment set-points, depending on the incident wind conditions, to match $\bm{\gamma^*}$, the model-optimal wake steering control strategy.
We tabulate $\bm{\gamma^*}$ for three independent input variables: wind speed, wind direction, and turbulence intensity.
Details regarding the calibration and optimization of the flow control model are in the {\bf Supplemental Information.}
During this experiment, we construct datasets for the wind farm power production in both wake steering and baseline yaw aligned control.
We switch between wake steering control and baseline control every $2.5$ hours (details in {\bf Methods}).
We performed the energy maximization experiment from December 2020 until March 2021.

\section*{Results}

\subsection*{Model validation with fixed yaw angles}

We consider the power production of the three turbine array for incident wind speeds in Region II of the power curve (Region II wind speeds are nominally between $4$ and $12$ m/s where turbines maximize their coefficients of power \cite{burton2011wind}).
The data filtering methodology is described in {\bf Methods.}
We normalize the power production for each turbine at each one-minute averaged instance by the power production of the adjacent reference turbine.
The model is calibrated using Eq. \ref{eq:cali} for instances in which turbine 1 commands zero yaw, $\gamma_1=0^\circ$.
Given the calibration based on standard control, the model is then used to predict the power production given the various yaw misalignment angles which were implemented for turbine 1 (Figure \ref{fig:setup}(c)).

We first consider northerly incident wind ($\alpha=0^\circ$) for which the three turbine array is slightly offset from alignment.
Both the measured and model predicted power for yaw misalignments of $\gamma_1=-20^\circ$, $\gamma_1=+20^\circ$, and yaw aligned operation $\gamma_1=0^\circ$ are shown in Figure \ref{fig:contours}.
For baseline operation ($\gamma_1=0^\circ$), the power of turbines $2$ and $3$ are approximately $40\%$ of the power of turbine $1$ (Figure \ref{fig:contours}(b)).

Our model predicts that negative yaw misalignment would reduce the array power production compared to baseline yaw aligned control as the wake of turbine $1$ will be steered towards the downwind turbines (Figure \ref{fig:contours}(a)).
Overall, a negative yaw misalignment of $\gamma_1=-20^\circ$ reduces sum of power production for the three turbines by $5.1\%$.
The power production for turbine $1$ decreases due to its yaw misalignment.
The power production for turbine $2$ slightly increases due to the reduction of the thrust of turbine $1$.
The power production for turbine $3$ is decreased for $\gamma_1=-20^\circ$, compared to $\gamma_1=0^\circ$.

The wake model velocity and power production predictions for $\gamma_1=+20^\circ$ and $\alpha=0^\circ$ are shown in Figure \ref{fig:contours}(c).
For $\gamma_1=+20^\circ$, the measured productions for the downwind turbines 2 and 3 are statistically significantly higher than for $\gamma_1=0^\circ$, as indicated by the $95\%$ confidence intervals.
Interestingly, the power for turbine $1$ is slightly higher with positive yaw misalignment, compared to yaw aligned operation, for both field measurements and model predictions.
This slight increase has been observed in other studies, depending on the incident wind conditions and turbine control system \cite{doekemeijer2021field, howland2020influence}.
Collectively, the yaw misalignment of turbine 1 to $\gamma_1=+20^\circ$ increased total power for the three turbine array by $28.6\%$.

The wake model predictions for all three yaw misalignment conditions shown in Figure \ref{fig:contours} are within or nearly within the $95\%$ confidence intervals of the field data, validating the predictive accuracy of the flow control model proposed in this study.
An outstanding question is whether the model can predict the true optimal yaw misalignment values for the utility-scale wind farm.
Here, we have yaw misaligned turbine $1$ for values between $\gamma_1=-25^\circ$ and $\gamma_1=25^\circ$ ($|\gamma_1|>25^\circ$ was not tested for loads limitations \cite{damiani2018assessment}).
The value of $\gamma_1$ that maximizes the array power depends on the wind direction $\alpha$.
We compare the values of $\gamma_1$ that maximize the array power in the utility-scale wind farm data to the values of $\gamma_1$ that maximize the power predictions from the flow control model for different inflow wind directions $\alpha$ (positive and negative $\alpha$ are northeasterly and northwesterly winds, respectively).
The results for $\alpha=-7.5^\circ, \ -5^\circ,$ and $0^\circ$ are shown in Figure \ref{fig:yaw_offset_1_2_3}(a-c) (all $6$ values of $\alpha$ with sufficient collected data are shown in the {\bf Supplemental Information}).
The value of $\gamma_1$ that maximizes power in the flow control model is within $5^\circ$ of the value that maximizes power in the field data for $5/6$ values of $\alpha$.
Based on the observed trend, it is possible that the optimal yaw in field conditions is $|\gamma|>25^\circ$ for $\alpha=-7.5^\circ$, but we do not have access to these data given the experimental setup and we do not expect $|\gamma^*|$ to be much greater than $30^\circ$ given the high power loss at turbine $1$ for large yaw \cite{howland2020influence}.
The results demonstrate that the flow control model is able to accurately estimate the true optimal yaw misalignment angles in utility-scale wind farms with sufficient precision to implement wake steering given typical yaw resolutions ($5^\circ$, \cite{howland2020influence}).

\begin{figure}
    \centering
    \includegraphics[width=0.75\linewidth]{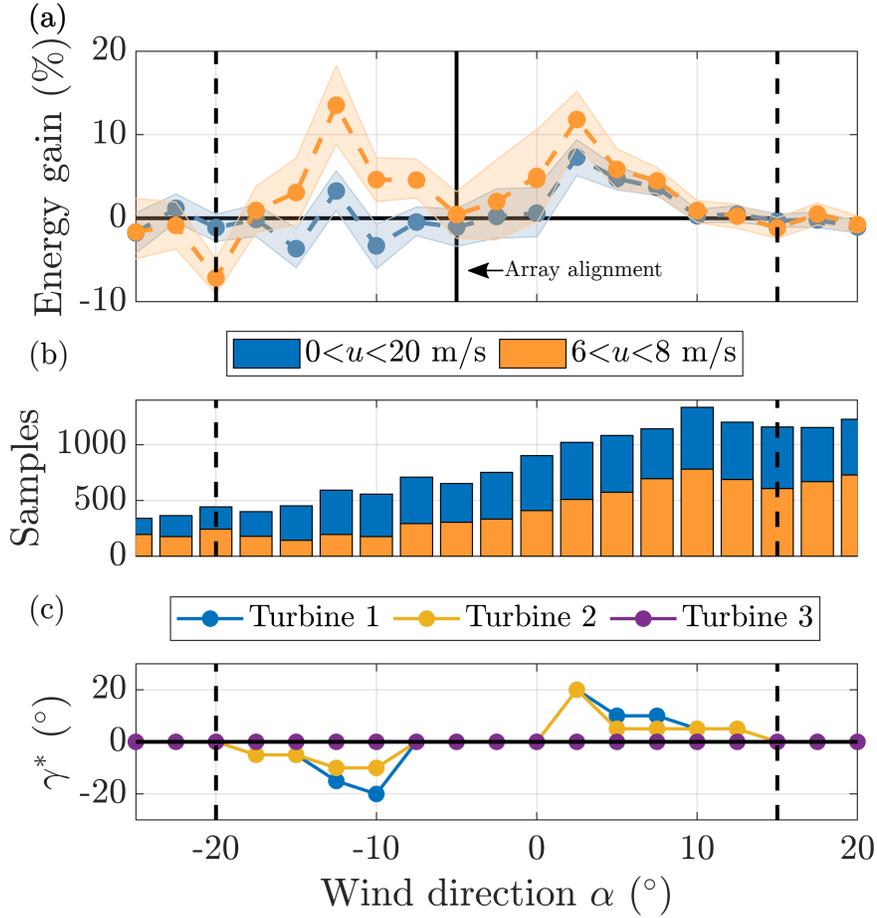}
    \caption{(a) Energy gain from the farm energy maximization wake steering field experiment. 
    The results for wind speeds between $6$\textless$u$\textless$8~\mathrm{m/s}$ and $0$\textless$u$\textless$20~\mathrm{m/s}$ (full wind speed range) are shown in orange and blue, respectively.
    The shaded region corresponds to $95\%$ confidence intervals around the mean estimated with bootstrapping.
    The wind turbines are approximately aligned for northwesterly inflow ($\alpha\approx-5^\circ$), shown with a solid black vertical line.
    (b) Number of unique one minute averaged data samples collected for each wind direction for wake steering operation.
    The number of data points collected for baseline yaw aligned operation is similar.
    (c) The yaw misalignment set-points applied to the farm which maximize wind farm power production in the flow control model $\gamma^*$ for each turbine.
    The optimal yaw angles $\gamma^*$ depend on the incident wind conditions of: wind speed, wind direction, and turbulence intensity.
    The yaw angles for incident wind speeds of $5\pm1~\mathrm{m/s}$ and turbulence intensities of $7.5\pm1.25~\%$ are shown.
    The full yaw misalignment lookup table is provided in the {\bf Supplemental Information}.}
    \label{fig:energy}
\end{figure}

\subsection*{Farm energy maximization}

We consider the impact of wake steering control on wind farm energy production.
The previous experiment demonstrated that wake steering can increase or decrease power production depending on the incident wind conditions and the prescribed yaw misalignment values (see Figure \ref{fig:yaw_offset_1_2_3}).
In this second experiment, we set the yaw misalignment values for each turbine in the farm to the model-optimal result $\bm{\gamma^*}$ depending on the incident wind speed, direction, and turbulence intensity.
Wake steering control is active ($\bm{\gamma}\neq0$) for $-17.5^\circ$\textless$\alpha$\textless$12.5^\circ$, since these are the wind directions with wake losses.
Although the commanded yaw is zero outside of $-17.5^\circ$\textless$\alpha$\textless$12.5^\circ$, our analysis considers $-20^\circ$\textless$\alpha$\textless$15^\circ$ to account for hysteresis in the yaw controller \cite{fleming2019initial}.

The effect of wake steering on the array power production will depend on the prescribed yaw misalignment and the incident wind conditions.
For example, we expect wake steering to have a relatively small effect on wind farm power production for high wind speeds since wake losses are relatively small.
However, higher speeds contribute more to wind farm energy than lower speeds.
To assess the collective benefit of wake steering, we assess the impact of wake steering on the array energy production.
The gain in energy from wake steering control, compared to baseline individual control, is computed using Eqs. \ref{eq-er}-\ref{eq:ratio-gain} and is described in {\bf Methods}.

First, we consider the gain in energy for each independent wind direction.
We first consider wind speeds between $6$ and $8$ m/s, the center of Region II of the power curve, selected to ensure we do not cross into other control regions.
These were the wind speeds of interest in the model validation experiment (see previous section).
The energy gain from wake steering control is shown for each wind direction in Figure \ref{fig:energy}(a).
The number of unique one-minute averaged data samples collected for wake steering operation for each $\alpha$ is shown in Figure \ref{fig:energy}(b) (similar numbers of samples were collected for aligned control).
Wake steering results in statistically significant gains in energy when active ($\bm{\gamma}\neq0$).

Considering the full wind speed range encountered ($0$ and $20$ m/s), the energy gain results are qualitatively similar although the magnitude of gain is decreased.
This is expected, since wake steering control will not be active for higher wind speeds given limited wake losses.
The energy gain from wake steering control decreased more for $\alpha<-5^\circ$ than for $\alpha>-5^\circ$.
This decrease is primarily associated with an energy decrease for the simultaneously waked and yaw misaligned turbine $2$ at moderate wind speeds between $8$ and $10$ m/s, as detailed further in the {\bf Supplemental Information.}

We compute the energy gain over a sector of wind directions using Eq. \ref{eq:sector_gain} ({\bf Methods}).
The energy gains over two wind direction sectors for wind speeds between $6$\textless$u$\textless$8$ m/s (Region II) are shown in Table \ref{tab:energy}.
For Region II operation, wake steering has a statistically significant increase for both the incident wind directions where wake steering is active ($2.7\%\pm0.8\%$) and the full wind direction distribution between $-180^\circ$ \textless $\alpha$ \textless $180^\circ$ ($1.0\%\pm0.6\%$).
Since wake steering is only active for wind directions which lead to wake interactions, the degree to which wake steering increases total energy production depends on the fraction of time in which the wind is oriented in these directions.
In this experiment, wake steering, which was only active between $-17.5^\circ$ and $12.5^\circ$, increased energy to a sufficient degree that it resulted in a statistically significant increase in the energy production in Region II for the full wind direction distribution ($-180^\circ$ to $180^\circ$).

\begin{table}
\small
\caption{Wind farm energy gain from wake steering over wind direction and wind speed sectors with $95\%$ confidence intervals from bootstrapping.
The wind farm energy gain is estimated using Eq. \ref{eq:sector_gain} ({\bf Methods}).} 
\centering 
\begin{tabular}{c|cc} 
Wind speed & \multicolumn{2}{c}{Incident wind direction sector} \\ 
(m/s) & -$20^\circ$ \textless $\alpha$ \textless $15^\circ$ & $-180^\circ$ \textless $\alpha$ \textless $180^\circ$ \\
\hline  
$6$\textless$u$\textless$8$
& 
${\color{seagreen} +2.7\%} \pm 0.8\%$
&
${\color{seagreen} +1.0\%} \pm 0.6\%$
\\
 $0$\textless$u$\textless$20$ (full $u$ range)
& 
${\color{seagreen} +1.0\%} \pm 0.5\%$
&
${\color{seagreen} +0.3\%} \pm 0.3\%$
\\
\hline 
\end{tabular}
\label{tab:energy}
\end{table}

The sector energy gains for the full wind speed range of $0$\textless$u$\textless$20$ m/s are shown in Table \ref{tab:energy}.
For the wind directions for which wake steering is active, the energy production increase due to wake steering is statistically significant ($1.0\%\pm0.5\%$).
For the full wind speed and direction range, the energy production gain in this experiment is $0.3\% \pm 0.3\%$.

\section*{Conclusions}

We demonstrate the potential for collective wind farm control at utility-scale.
In this study, we develop a predictive physics-based, data-assisted wind farm flow control model.
To validate the model, we design a multi-month field experiment at a utility-scale wind farm.
We demonstrate that the optimization-oriented flow control model is able to predict the optimal yaw misalignment angles for the utility-scale wind turbine array within $\pm 5^\circ$ for most conditions.

We design a wake steering strategy by maximizing the power production of the wind turbine array depending on the incident atmospheric conditions in the validated model.
For wind directions for which wake steering control is active and Region II wind speeds, the array energy production is increased by $2.7\% \pm 0.8\%$.
For the same wind directions, considering the full wind speed range (including when wake steering is not active at high winds), the energy gain is $1.0\% \pm 0.5\%$.
For Region II wind speeds, the energy gain across all incident wind directions is $1.0\% \pm 0.6\%$.
Finally, for all wind speeds and directions, the energy gain is $0.3\% \pm 0.3\%$.

The gain in energy production at a commercial wind farm depends on the farm layout and the site-specific wind conditions.
Given the incident winds observed during this experiment at the Indian farm, the energy gain is $1.0\%$ and $0.3\%$ for Region II and all wind speeds respectively.
For wind generation to produce a larger fraction of energy worldwide, wind farm development will have to expand to onshore regions of lower mean wind speed and offshore environments, where wake losses are anticipated to be higher \cite{bodini2019us}.
In this study, we demonstrate that wake steering control is a mechanism by which wind farm energy can be increased.
We also demonstrate that the proposed flow control modeling framework is able to predict the gain in production from wake steering, enabling its use in wind farm control and potentially farm design \cite{fleming2016wind}.

Further improvements in active wake control are necessary to increase the efficacy of wake steering, such as closed-loop control \cite{howland2021optimal,doekemeijer2020closed,howland2020optimal,ciri2017model} or machine learning based methods \cite{stanfel2021proof}.
Additionally, further improvements in flow control models are necessary to capture three dimensional effects such as wake curling \cite{howland2016wake, martinez2021curled, bastankhah2021analytical}, veer \cite{howland2020influence, abkar2018analytical}, and stability \cite{howland2021optimal, pena2014modeling}.

To achieve climate goals, renewables must produce most of the demanded global energy.
Wake steering is especially effective at increasing farm production at low wind speeds.
Since annual wind speeds in India, and other emerging economies, are relatively low \cite{purohit2009wind}, wake steering is a promising mechanism to increase renewable energy production.
Beyond higher annual energy production, low-carbon grids will require more predictable and controllable wind energy with reduced intermittency.
Wake steering contributes to our ability to actively control the production from wind farms.
Improvements in our ability to model, control and design wind energy generation will facilitate a more rapid energy grid decarbonization.

\section*{Methods} \label{sec: Methods}

\subsection*{Predictive farm flow control model}

With the goal of selecting the optimal operational strategy for utility-scale wind turbines, we develop a wind farm flow control model which accounts for the variable atmospheric conditions.
To remain tractable for  controls-oriented modeling, the wind farm flow model must be computationally efficient \cite{boersma2016control}.
But the model must be sufficiently accurate to ensure that the model-optimal yaw angles are closely representing the true optimal yaw angles in complex field conditions.

In wind farms, there are limited sensors which measure the freestream wind flow at multiple heights (see Figure \ref{fig:model}).
Sensors which measure the wind flow up to the wind turbine hub height, which is approximately $100$ meters in altitude, are expensive \cite{eberle2019nrel}.
Such sensors tend to be placed sparsely and suboptimally in wind farms \cite{annoni2018sparse} and are only useful for measuring the incident wind profiles, rather than velocities in turbine wakes.
In summary, we generally have access to the incident wind profile through LiDAR or meteorological tower measurements, but we require a wake model to predict the wind flow field in the wake of turbines as a function of their operational strategy (see Figure \ref{fig:model}).
Further discussion of the role of the atmospheric conditions on the wind farm power is given in the {\bf Supplemental Information.}

We model the power production of yaw misaligned wind turbines given the incident wind conditions using blade element theory.
We denote the modeled power of upwind (unwaked) turbines as $\hat{P}_u$.
Previous approaches have used a tuned empirical model to estimate the power production of yaw misaligned turbines, where the power of the yaw misaligned turbine is estimated as $\hat{P}_u(\gamma)=\hat{P}_u(\gamma=0)\cdot\cos^{P_p}(\gamma)$, where $P_p$ is a tuned parameter \cite{gebraad2016wind}. 
This empirical approach leads to significant errors depending on the incident winds in the ABL \cite{howland2020influence}.
Given variations in the wind speed and direction as a function of height $z$ in the ABL, we calculate the forces on the blade at each location in the wind turbine rotor area as a function of the radial $r$ and azimuthal $\theta$ positions.
The wind turbine power production is modeled as the integration of the incremental torque over the rotor area multiplied by the angular velocity
\begin{equation}
\hat{P}_u = \Omega \int_0^{2\pi} \int_0^r r df_\tau(r,\theta), 
\end{equation}
where $\Omega$ is the blade angular velocity and $df_\tau$ is the tangential force at a particular blade section.
The tangential forces depend on characteristics of the wind turbine, including aerodynamic properties of the blades, the imposed yaw misalignment, and the incident wind profiles \cite{howland2020influence}.
More details are provided in the {\bf Supplemental Information.}

We model the power production of downwind, waked turbines ($\bm{\hat{P}_w}$) using a computationally efficient analytical wake model developed by Howland \textit{et al.} \cite{howland2019wind}, which is based on a lifting line wake deflection model \cite{shapiro2018modelling}.
The full wake model is described in the {\bf Supplemental Information}.
We leverage modified linear wake superposition and an analytical secondary steering model \cite{howland2021influence}, which models the effect of the nonzero lateral velocity at turbines downwind of a yawed turbine \cite{king2020controls}.
To predict wind turbine power production in a computationally efficient fashion, analytical wake models parameterize wake and ABL turbulence with an unknown wake spreading rate $\bm{k_w}$ \cite{stevens2017flow, howland2021optimal, meneveau2019big}.
Traditional approaches have tuned the wake model parameter to idealized data from simulations \cite{niayifar2016analytical} which neglects the site-specific effects on the model parameters \cite{howland2019wind} and the effects of the stability of the ABL \cite{howland2021optimal}.
Instead, we estimate the wake model parameters using optimization-based inverse problem techniques \cite{stuart2010inverse}, where we minimize the mean squared difference between the model predictions and field calibration data
\begin{equation}
\bm{k_w}^* = \underset{\bm{k_w}}\argmin~\overline{(\bm{\hat{P}_w}(\bm{k_w}, \bm{\gamma=0}) - \bm{P_w( \bm{\gamma}=0)})^2}.
\label{eq:cali}
\end{equation}
We calibrate the wake spreading rate using an offline variant of the ensemble Kalman filter \cite{howland2021optimal, howland2020optimal, evensen2003ensemble}.
We calibrate the wake spreading rate using site- and wind condition-specific historical data collected where the wind farm is in standard, individual operation ($\bm{\gamma=0}$).
More details regarding the calibration are provided in the {\bf Supplemental Information.}
The calibrated model parameters are then used to predict the power production for the wind farm in wake steering configurations ($\bm{\gamma\neq0}$).
The optimal wind farm control strategy with $N_t$ turbines is predicted using the flow control model as
\begin{equation}
\bm{\gamma}^* = \underset{\bm{\gamma}}\argmax \sum_{i=1}^{N_t} \hat{P}_i(\bm{\gamma}).
\label{eq:opti}
\end{equation}
We estimate $\bm{\gamma}^*$ using gradient-based optimization of the flow control model \cite{howland2019wind}.
More details regarding the optimization are provided in the {\bf Supplemental Information.}

\subsection*{Fixed yaw misalignment experiment data filtering}

For the fixed yaw misalignment experiment, turbine $1$ (Figure \ref{fig:setup}(b)) implements a yaw misalignment set-point sequence between $\gamma = -25^\circ$ and $+25^\circ$ (Figure \ref{fig:setup}(c)).
The yaw misalignment set-points are active for Region II wind speeds, between approximately $4$ m/s and $12$ m/s.
Since the wind direction in the ABL is constantly evolving, the standard wind turbine yaw control system is used to track the yaw misalignment set-point \cite{burton2011wind}.
The yaw set-point tracking ability of the turbine model of interest was validated by Howland \textit{et al.} \cite{howland2020influence}.
Data are aggregated in $1~\mathrm{min}$ averages.

For the fixed yaw misalignment experiment, we consider conditional averages of the power ratio for the turbines of interest.
We characterize the freestream wind conditions using measurements made by the reference turbine (Figure \ref{fig:setup}(b)).
We focus on wind speeds between $6$ and $8$ m/s.
Since the wake interactions depend on the turbulence intensity \cite{howland2021optimal}, we restrict the analysis to turbulence intensities between $2.5\%$ and $7.5\%$.
We characterize the realized yaw misalignment of turbine $1$ based on the nacelle-mounted measurements of the relative wind direction of turbine $1$, measured by a sonic anemometer.
We consider yaw misalignment for values between $-25^\circ$ and $25^\circ$ with a step size of $5^\circ$.
The wind direction is considered for values between $-10^\circ$ to $2.5^\circ$ with a step size of $2.5^\circ$.
The results for all incident wind directions are shown in the {\bf Supplemental Information}.

\subsection*{Energy gain analysis methodology}

To estimate the impact of the wake steering control on energy production, we compute the energy ratio, which is the ratio of the energy production of the test turbine of interest to the energy production of the reference turbine.
The energy ratio is given by
\begin{equation}
E_r = \frac{\sum_{i=1}^{N_b} w_i \overline{P}_i^{\mathrm{test}}}{ \sum_{i=1}^{N_b} w_i \overline{P}_i^{\mathrm{ref}}},
\label{eq-er}
\end{equation}
where $N_b$ is the number of wind speed bins, $w_i$ is the wind condition specific weighting factor, and $\overline{P}$ is the mean power production within the particular wind speed bin.

The energy ratio is computed for each wind direction.
Additionally, the energy ratio is computed for both baseline yaw aligned control ($\mathrm{base}$) and for wake steering control ($\mathrm{steer}$).
In the field experiment of Fleming {\it et al.} (2019) \cite{fleming2019initial}, the weighting factor was defined as the probability distribution over the wind speeds for a particular wind direction for the alternative control method.
The weighting factors are
\begin{align}
w_i^{\mathrm{base}} &= \frac{n_i^{\mathrm{steer}}}{N^{\mathrm{steer}}} &  
w_i^{\mathrm{steer}} &= \frac{n_i^{\mathrm{base}}}{N^{\mathrm{base}}},
\label{eq:weights}
\end{align}
where $w_i^{\mathrm{base}}$ and $w_i^{\mathrm{steer}}$ are the weights ($w_i$) used in the energy ratio calculation for the wake steering and the baseline yaw aligned cases, respectively. 
The number of samples in the wind speed bin for a particular wind direction is given by $n_i$.
The number of total samples for a particular wind direction is given by $N$.
Note that the weights are switched to the opposite experimental window, as discussed in Fleming et al. (2019) \cite{fleming2019initial}, to alleviate the relatively large wake steering energy gains at low wind speeds. 
In the limit of identical wind conditions in wake steering and yaw aligned control, the effect of switching the weights becomes null. 
We test the influence of different weighting approaches in the {\bf Supplemental Information}.
We find that the different weighting methods do not have a statistically significant impact on our energy gain results.
The fractional gain in energy for each wind direction is
\begin{equation}
R = \frac{E_{r}^{\mathrm{steer}}}{E_{r}^{\mathrm{base}}}.
\label{eq:ratio-gain}
\end{equation}

The influence of wake steering on the energy production over a sector of wind directions can be approximated using the energy ratio. 
The sector energy gain is
\begin{eqnarray}
\mathrm{G}_r = \sum_{i=1}^{N_\alpha} W_i^{\mathrm{ref}} R_i,
\label{eq:sector_gain}
\end{eqnarray}
where $N_\alpha$ is the number of wind directions considered in the sector energy analysis and $R_i$ is the energy ratio gain for each wind direction (Eq. \ref{eq:ratio-gain}).
The weights are
\begin{equation}
W_i^{\mathrm{ref}} = \frac{\overline{E}_i^{\mathrm{ref}} N_i^{\mathrm{ref}}}{\sum_{i=1}^{N_\alpha} \overline{E}_i^{\mathrm{ref}} N_i^{\mathrm{ref}}},
\end{equation}
which represents the fraction of energy produced in each wind direction (accounting for the wind rose).
The averaged energy produced by the reference turbine for each wind direction $\overline{E}_i^{\mathrm{ref}}$ and $N_i^{ref}$ is the number of samples for each wind direction.
Confidence intervals for each statistic are calculated using bootstrapping.

\subsection*{Farm energy maximization data filtering}

For the wake steering farm energy maximization experiment, we switch between wake steering control and baseline yaw aligned control with a period of $2.5~\mathrm{h}$.
We selected a switching period of $2.5~\mathrm{h}$ to balance the aim of collecting data in similar conditions while limiting the number of transitions between control strategies.
When baseline yaw aligned control is activated, all turbines implement a strategy of $\gamma=0^\circ$.
When wake steering is activated, turbines $1$ and $2$ implement a dynamic yaw misalignment sequence based on the incident wind conditions (Figure \ref{fig:energy}(c)) and turbine $3$ implements a yaw aligned strategy of $\gamma=0^\circ$.
As with the fixed yaw experiment, the standard wind turbine yaw control system is used to track the intended yaw misalignment set-point, which now varies as a function of the incident wind conditions.
Data are aggregated in $1~\mathrm{min}$ averages.

We compute the energy ratio (Eq. \eqref{eq-er}) using the data for the two different control approaches (wake steering and baseline yaw aligned).
We consider wind speed and direction ranges as specified by Table~\ref{tab:energy} with a step size of $1~\mathrm{m/s}$.
We consider only wind directions with more than $20$ data points in both control cases.
We consider only wind speed subsets of a wind direction bin with more than $5$ data points in both control cases.
Since high turbulence was infrequent in these experiments, gives high noise-to-signal ratio and turbine loads \cite{damiani2018assessment}, and wake losses are minimal in inflow with high turbulence, turbulence intensities less than $20\%$ are considered for wake steering.
We only consider data in which all four turbines (three test turbines and reference turbine) are active and operating normally.

\FloatBarrier


\bibliography{scibib}
\bibliographystyle{unsrt} 


\end{document}